# Sufficient and necessary conditions for generalized Hölder's inequality in p-summable sequence spaces

**Al A Masta*[1], S Fatimah[2], A Arsisari[3], Y Y Putra[4] and F Apriani[5]**

[1,2]Department of Mathematics Education, Universitas Pendidikan Indonesia, Jl. Dr. Setiabudi, Bandung, Indonesia

[3,4,5]Department of Mathematics Education, Sekolah Tinggi Keguruan dan Ilmu Pendidikan Muhammadiyah Bangka Belitung, Jl. KH. Ahmad Dahlan, Kabupaten Bangka Tengah, Bangka Belitung, Indonesia

*Corresponding author's e-mail: alazhari.masta@upi.edu

**Abstract**. In this paper, we discuss the generalized Hölder's inequality in p-summable sequence spaces. In particular, we shall prove sufficient and necessary conditions for generalized Hölder's inequality in those spaces. One of the keys to prove our results is to estimate the norms of characteristic sequences in $\mathbb{R}$.

## 1. Introduction

In mathematics, Lebesgue spaces are one of important topics, particularly in real and functional analysis. There are two kinds of Lebesgue spaces which are 'continuous' Lebesgue spaces denoted by $L_p$ and p-summable sequence spaces denoted by $\ell_p$. Many researchers have studied Lebesgue spaces and its generalization over few decades (see [1 - 12], etc.). For example, in 2016, Masta, *et al.* [2] presented the sufficient and necessary conditions for generalized Hölder's inequality in 'continuous' Lebesgue spaces and in 'continuous' Orlicz – Morrey spaces. Recently, Ifronika, *et al.* [3] also obtained the sufficient and necessary conditions for generalized Hölder's inequality in 'continuous' Morrey spaces, in 'continuous' generalized Morrey spaces, and in their weak type.

Motivated by these results, we are interested in discussing the p-summable sequence spaces. In particular, we will prove the sufficient and necessary conditions for Hölder's inequality in these spaces.

First, we recall the definition of p-summable sequence spaces. Let $1 \leq p < \infty$, the p-summable sequence $\ell_p(\mathbb{R})$ is the set of sequences $X = (x_n) \subseteq \mathbb{R}$ such that

$$\|X\|_{\ell_p(\mathbb{R})} := \left(\sum_{n=1}^{\infty} |x_n|^p\right)^{\frac{1}{p}} < \infty.$$

Note that, $\ell_p(\mathbb{R})$ is a Banach space with respect to the norm $\|\cdot\|_{\ell_p(\mathbb{R})}$. Next, for $p = \infty$, $\ell_\infty(\mathbb{R})$ is defined as the set of sequences $X = (x_n) \subseteq \mathbb{R}$ such that

$$\|X\|_{\ell_\infty(\mathbb{R})} := \max_{n\in\mathbb{N}}|x_n| < \infty.$$

Notice that, $\ell_\infty(\mathbb{R})$ is also Banach spaces equipped with the norm $\|\cdot\|_{\ell_\infty(\mathbb{R})}$.

The rest of this paper is organized as follows. In Section 2, we presented some lemmas which useful for obtain our results. The main results are presented in Section 3. In Section 3, we state the sufficient and necessary conditions for generalized Hölder's inequality in p-summable sequence spaces.

## 2. Methods

To obtain the sufficient and necessary conditions for generalized Hölder's inequality in p-summable sequence spaces, we use the norms of the characteristic sequences in $\mathbb{R}$ and some lemmas as in the following.

**Lemma 1.** [4] *Let* $m \in \mathbb{Z}$ *and* $N \in \{0,1,2,3,\ldots\}$, *write* $S_{m,N} := \{m-N, \ldots m, \ldots, m+N\}$. *Let*

$$\xi_k^{m,N} := \begin{cases} 1, & \text{if } k \in S_{m,N} \\ 0, & \text{otherwise} \end{cases},$$

*then there exists* $C > 0$ *(independent of $m$ and $N$) such that*

$$(2N+1)^{1/p} \leq \|\xi_k^{m,N}\|_{\ell_p(\mathbb{R})} \leq C(2N+1)^{1/p}$$

*for every* $N \in \{0,1,2,3,\ldots\}$.

**Lemma 2.** [13] *Let $x_i$ be positive real numbers for $i = 1,2,3,\ldots,m$. If $1 \leq p_1, p_2, \ldots, p_m, p < \infty$ satisfy the condition $\sum_{i=1}^{m} \frac{1}{p_i} = \frac{1}{p}$, then wen have*

$$\left(\prod_{i=1}^{m} x_i\right)^p \leq \sum_{i=1}^{m} \frac{p}{p_i} x_i^{p_i}.$$

**Corollary 3.** *Let* $1 \leq p_1, p_2, p < \infty$ *satisfy the condition* $\frac{1}{p_1} + \frac{1}{p_2} = \frac{1}{p}$, $X = (x_n) \in \ell_{p_1}(\mathbb{R})$ *and* $Y = (y_n) \in \ell_{p_2}(\mathbb{R})$. *If* $\sum_{n=1}^{\infty}|x_n|^{p_1} = \sum_{n=1}^{\infty}|y_n|^{p_2} = 1$, *then* $\sum_{n=1}^{\infty}|x_n y_n|^p \leq 1$.

Proof.
Let $\frac{1}{p_1} + \frac{1}{p_2} = \frac{1}{p}$, $X = (x_n) \in \ell_{p_1}(\mathbb{R})$ and $Y = (y_n) \in \ell_{p_2}(\mathbb{R})$. By using Lemma 2, we have

$$\sum_{n=1}^{\infty}|x_n y_n|^p \leq \frac{p}{p_1}\sum_{n=1}^{\infty}|x_n|^{p_1} + \frac{p}{p_2}\sum_{n=1}^{\infty}|y_n|^{p_2} = \frac{p}{p_1} + \frac{p}{p_2} = 1. \blacksquare$$

## 3. Results and Discussion

First, we present sufficient and necessary conditions for Hölder's inequality in $\ell_p(\mathbb{R})$ space in the following theorem.

**Theorem 4.** *Let* $1 \leq p, p_1, p_2 < \infty$.
(1) *If* $\frac{1}{p_1} + \frac{1}{p_2} = \frac{1}{p}$, *then* $\|XY\|_{\ell_p(\mathbb{R})} \leq \|X\|_{\ell_{p_1}(\mathbb{R})} \|Y\|_{\ell_{p_2}(\mathbb{R})}$ *for every* $X \in \ell_{p_1}(\mathbb{R})$ *and* $Y \in \ell_{p_2}(\mathbb{R})$.
(2) *If* $\|XY\|_{\ell_p(\mathbb{R})} \leq \|X\|_{\ell_{p_1}(\mathbb{R})} \|Y\|_{\ell_{p_2}(\mathbb{R})}$, *for every* $X \in \ell_{p_1}(\mathbb{R})$ *and* $Y \in \ell_{p_2}(\mathbb{R})$, *then* $\frac{1}{p_1} + \frac{1}{p_2} \geq \frac{1}{p}$.

Proof.

(1). Let $\frac{1}{p_1}+\frac{1}{p_2}=\frac{1}{p}$, $X=(x_n)\in \ell_{p_1}(\mathbb{R})$ and $Y=(y_n)\in \ell_{p_2}(\mathbb{R})$. First, suppose that $\sum_{n=1}^{\infty}|x_n|^{p_1}=A$ and $\sum_{n=1}^{\infty}|y_n|^{p_2}=B$. By setting $x'_n=\frac{x_n}{A^{1/p_1}}$ and $y'_n=\frac{y_n}{B^{1/p_2}}$, we have $\sum_{n=1}^{\infty}|x'_n|^{p_1}=\frac{\sum_{n=1}^{\infty}|x_n|^{p_1}}{A}=1$ and $\sum_{n=1}^{\infty}|y'_n|^{p_2}=\frac{\sum_{n=1}^{\infty}|y_n|^{p_2}}{B}=1$. By using Corollary 3, we have

$$1\geq \sum_{n=1}^{\infty}|x'_n y'_n|^p = \sum_{n=1}^{\infty}\left|\frac{x_n y_n}{A^{\frac{1}{p_1}}B^{\frac{1}{p_2}}}\right|^p = \frac{1}{A^{\frac{p}{p_1}}B^{\frac{p}{p_2}}}\sum_{n=1}^{\infty}|x_n y_n|^p.$$

Since $1\geq \frac{1}{A^{\frac{p}{p_1}}B^{\frac{p}{p_2}}}\sum_{n=1}^{\infty}|x_n y_n|^p$ is equivalent to $\sum_{n=1}^{\infty}|x_n y_n|^p \leq A^{\frac{p}{p_1}}B^{\frac{p}{p_2}}$ we obtain

$$\left(\sum_{n=1}^{\infty}|x_n y_n|^p\right)^{1/p} \leq \left(\sum_{n=1}^{\infty}|x_n|^{p_1}\right)^{\frac{1}{p_1}}\left(\sum_{n=1}^{\infty}|y_n|^{p_2}\right)^{\frac{1}{p_2}}.$$

So, we have $\|XY\|_{\ell_p(\mathbb{R})} \leq \|X\|_{\ell_{p_1}(\mathbb{R})}\|Y\|_{\ell_{p_2}(\mathbb{R})}$.

(2). Now, assume that $\|XY\|_{\ell_p(\mathbb{R})} \leq \|X\|_{\ell_{p_1}(\mathbb{R})}\|Y\|_{\ell_{p_2}(\mathbb{R})}$ holds, for every $X\in \ell_{p_1}(\mathbb{R})$ and $Y\in \ell_{p_2}(\mathbb{R})$. Take $X=Y=\xi_k^{m,N}$, by using Lemma 1, we have

$$(2N+1)^{1/p} \leq \|\xi_k^{m,N}\|_{\ell_p(\mathbb{R})} \leq \|\xi_k^{m,N}\|_{\ell_{p_1}(\mathbb{R})}\|\xi_k^{m,N}\|_{\ell_{p_2}(\mathbb{R})} \leq C^2(2N+1)^{\frac{1}{p_1}+\frac{1}{p_2}}$$

or $(2N+1)^{\frac{1}{p}-\left(\frac{1}{p_1}+\frac{1}{p_2}\right)} \leq C^2$ for every $N\in\{0,1,2,3,...\}$. Hence, we can conclude that $\frac{1}{p}\leq \frac{1}{p_1}+\frac{1}{p_2}$. ∎

**Theorem 5.** Let $1\leq p_1, p_2, p_3, ..., p_m, p < \infty$.
(1) If $\sum_{i=1}^{m}\frac{1}{p_i}=\frac{1}{p}$, then $\|\prod_{i=1}^{m}X_i\|_{\ell_p(\mathbb{R})} \leq \prod_{i=1}^{m}\|X_i\|_{\ell_{p_i}(\mathbb{R})}$, for every $X_i \in \ell_{p_i}(\mathbb{R})$.
(2) If $\|\prod_{i=1}^{m}X_i\|_{\ell_p(\mathbb{R})} \leq \prod_{i=1}^{m}\|X_i\|_{\ell_{p_i}(\mathbb{R})}$, for every $X_i \in \ell_{p_i}(\mathbb{R})$, then $\sum_{i=1}^{m}\frac{1}{p_i}\geq \frac{1}{p}$.

Proof.
(1). Let $\sum_{i=1}^{m}\frac{1}{p_i}=\frac{1}{p}$, $X_i=(x_{n,i})\in \ell_{p_i}(\mathbb{R})$ for $i=1,2,3,...,m$. First, suppose that $\sum_{n=1}^{\infty}|x_{n,i}|^{p_i}=A_i$ for every $i=1,2,...,m$. By setting $x'_{n,i}=\frac{x_{n,i}}{A_i^{1/p_i}}$, we have $\sum_{n=1}^{\infty}|x'_{n,i}|^{p_i}=\frac{\sum_{n=1}^{\infty}|x_{n,i}|^{p_i}}{A_i}=1$. By using Lemma 4, we have

$$\sum_{n=1}^{\infty}\left|\prod_{i=1}^{m}x_{n,i}\right|^p \leq \frac{p}{p_1}\sum_{n=1}^{\infty}|x_{n,1}|^{p_1} + \frac{p}{p_2}\sum_{n=1}^{\infty}|x_{n,2}|^{p_2} + \cdots + \frac{p}{p_m}\sum_{n=1}^{\infty}|x_{n,m}|^{p_m}$$

$$=\sum_{i=1}^{m}\frac{p}{p_i}=1.$$

So we have,

$$1 \geq \sum_{n=1}^{\infty}\left|\prod_{i=1}^{m} x'_{n,i}\right|^p = \sum_{n=1}^{\infty}\left|\prod_{i=1}^{m} \frac{x_{n,i}}{A_i^{\frac{1}{p_i}}}\right|^p = \frac{1}{\prod_{i=1}^{m} A_i^p} \sum_{n=1}^{\infty}\left|\prod_{i=1}^{m} x_{n,i}\right|^p.$$

Since $1 \geq \prod_{i=1}^{m} A_i^{-\frac{p}{p_i}} \left(\sum_{n=1}^{\infty}|\prod_{i=1}^{m} x_{n,i}|^p\right)$ is equivalent to $\sum_{n=1}^{\infty}|\prod_{i=1}^{m} x_{n,i}|^p \leq \prod_{i=1}^{m} A_i^{\frac{p}{p_i}}$, we obtain

$$\left(\sum_{n=1}^{\infty}\left|\prod_{i=1}^{m} x_{n,i}\right|^p\right)^{1/p} \leq \prod_{i=1}^{m}\left(\sum_{n=1}^{\infty}|x_{n,i}|^{p_i}\right)^{\frac{1}{p_i}}.$$

So, we have

$$\left\|\prod_{i=1}^{m} X_i\right\|_{\ell_p(\mathbb{R})} \leq C \left\|\prod_{i=1}^{m} X_i\right\|_{\ell_{p^*}(\mathbb{R})} \leq C \prod_{i=1}^{m} \|X_i\|_{\ell_{p_i}(\mathbb{R})}.$$

(2). Now, assume that $\|\prod_{i=1}^{m} X_i\|_{\ell_p(\mathbb{R})} \leq \prod_{i=1}^{m}\|X_i\|_{\ell_{p_i}(\mathbb{R})}$ holds for every $X_i \in \ell_{p_i}(\mathbb{R})$. Take $X_i = \xi_k^{m,N}$ for every $i = 1,2,3,\dots,m$, by using Lemma 1, we have

$$(2N+1)^{1/p} \leq \left\|\xi_k^{m,N}\right\|_{\ell_p(\mathbb{R})} \leq \prod_{i=1}^{m}\left\|\xi_k^{m,N}\right\|_{\ell_{p_i}(\mathbb{R})} \leq C^m (2N+1)^{\sum_{i=1}^{m}\frac{1}{p_i}}$$

or $(2N+1)^{\frac{1}{p} - \left(\sum_{i=1}^{m}\frac{1}{p_i}\right)} \leq C^m$ for every $N \in \{0,1,2,\dots\}$. Hence, we can conclude that $\frac{1}{p} \leq \sum_{i=1}^{m}\frac{1}{p_i}$. ∎

## 4. Conclusion

We have shown the sufficient and necessary conditions for generalized Hölder's inequality in $\ell_p(\mathbb{R})$ space, we can state that the condition $\frac{1}{p} \leq \sum_{i=1}^{m}\frac{1}{p_i}$ is a necessary conditions for generalized Hölder's inequality in $\ell_p(\mathbb{R})$ space.

**Acknowledgement.** The first and second authors is supported by Hibah Penguatan Kompetensi UPI 2018.